\documentclass{article}
\usepackage{mathrsfs}
\usepackage{graphicx}
\usepackage{amsmath}
\usepackage{amssymb}
 \usepackage{txfonts}

\newtheorem{thm}{Theorem}[section]

\newtheorem{lem}[thm]{Lemma}

\newtheorem{defi}[thm]{Definition}

\begin{document}

\begin{center}

{\Large \bf  On $c$-Normal Maximal Chains of Sylow Subgroups in Finite Groups\footnote{* Corresponding author: W. Meng (mlwhappyhappy@163.com).
W. Meng is supported by National Natural Science Foundation of China (12661001), Guangxi Natural Science
				Foundation Program (2026GXNSFAA00641090) and Guangxi Colleges and Universities Key Laboratory of Data Analysis and Computation.
}
}

\end{center}

\vskip0.5cm

\begin{center}

 Wei Meng*,  Feiyu Geng 

  a. School of Mathematics and Computing Science, Guilin University of Electronic Technology, Guilin, Guangxi, 541002, P.R. China.

b. Center for Applied Mathematics of Guangxi(GUET), Guilin, Guangxi, 541002, P.R. China.

  E-mails: mlwhappyhappy@163.com, 1300292725@qq.com.

\end{center}

\begin{abstract}
Let $G$ be finite group. A subgroup $H$ of $G$ is said to be  $c$-normal in
  $G$ if there exists a normal subgroup $N$ of $G$ such that $G = HN$ and $H\cap N\leq H_G$, where $H_G$ is the core of $H$ in $G$.
Let $P$ be a nontrivial $p$-group. A chain of subgroups $1=P_0 <P_1<\cdots P_{n-1}<P_n=P$   is called a maximal chain of $P$ provided that $|P_i:P_{i-1}|=p$ for $i=1,2,\ldots n-1$.  In this paper, we determine the structure of finite groups having $c$-normal subgroups of
 of a maximal chain of Sylow subgroups. We
obtain new characterizations of finite $p$-nilpotent and supersolvable groups
\end{abstract}




\section{Introduction}
Throughout this paper, all groups are finite. Let $G$ be a finite group and let $N$ be a normal subgroup of $G$. Then $NK=KN$ 
 for every subgroup $K$ of $G$.  Motivated by this basic commutativity property, many authors have generalized the notion of normality. For example, a subgroup $H$ of $G$ is called to be quasinormal in $G$ if  $HK=KH$ for every subgroup $K$ of $G$. 
 In this direction, Kegel \cite {KE} introduced the notion of $S$-quasinormal subgroups: a subgroup $H$ of a group $G$ is called $S$-quasinormal (also termed $S$-permutable, or $\pi$-quasinormal in some literature) in $G$ if $PH=HP$ for  every Sylow subgroup $P$ of $G$. This notion was later refined by Ballester-Bolinches and Pedraza-Aguilera \cite {BB}, who introduced $S$-quasinormally embedded subgroups. A subgroup
$H$ of $G$ is said to be $S$-quasinormally embedded (or $\pi$-quasinormally embedded) in $G$ if, for every prime $p\in\pi(H)$,  some Sylow
$p$-subgroup of $H$ is also a Sylow $p$-subgroup of an $S$-quasinormal subgroup of $G$. In a further generalization, Li et al. \cite {LI3} extended
$S$-permutability to $SS$-quasinormality. A subgroup $H$ of $G$ is called $SS$-quasinormal in $G$ if  there exists a subgroup $B$ of $G$  such that
$G=HB$ and $H$ permutes with every Sylow subgroup of $B$, i.e., $HP=PH$ for every Sylow subgroup $P$ of $B$. 
 The above generalized normality notions have a significant impact on the structure of finite groups. In particular, many authors have investigated the structure of finite groups under the assumption that certain subgroups of prime power order of $G$, or certain maximal subgroups of Sylow subgroups of $G$, possess these generalized normality properties in $G$,  and have thereby obtained numerous results (see \cite {BB,LI3,LI4, LI5,MENG1, MENG,QI2}).

Let $P$ be a $p$-subgroup of $G$.  A subgroup series  
\begin{center}
$\Omega:$ $1=P_0<P_1 <\cdots<P_i<\cdots<P_{n-1}< P_n=P$
\end{center}
is called a maximal chain  of $P$ if $P_{i-1}$ is a maximal subgroup of $P_i$ (i.e., $|P_i:P_{i-1}|=p$)  
for every $i \in\{1,2,\ldots,n\}$. 

Asaad \cite {ASA} established new characterizations of $p$-nilpotent and supersolvable finite groups by requiring each member of a maximal chain of Sylow $p$-subgroups to have a certain generalized normality property. Xu and Chen \cite {XU1} characterized
$p$-nilpotent groups and supersolvable groups in terms of the $\pi$-quasinormal embeddedness of a maximal chain of Sylow subgroups. Gao et al. \cite {GAO} also derived a characterization of $p$-nilpotent groups via the $\mathfrak{Z}$-permutability of the $p$-Sylowizers associated with a maximal chain of Sylow subgroups. More recently, Xu \cite {XU2} gave new characterizations of $p$-nilpotent and supersolvable groups by assuming that every member of a maximal chain of Sylow $p$-subgroups is $SS$-quasinormal.\\

Recall that a  subgroup $H$  of $G$ is said to be $c$-normal in
 $G$ if there exists a normal subgroup $N$ of $G$ such that $G=HN$ and $H\cap N\leq H_G$, where $H_G$ is the core of $H$ in $G$. The notion was introduced by Wang \cite{WANG}, and the influence of $c$-normal subgroups on the structure of finite groups has since been studied extensively(see \cite {GUO1, LI1, WEI1, WEI2}).

 Inspired by the above studies, we shall focus on finite groups $G$ in which every member of a maximal chain of Sylow subgroups is $c$-normal in $G$. For convenience,  we give the following definition.
 
 \begin{defi}
Let $G$ be a finite group and $P$ a $p$-subgroup of $G$. A maximal chain 
\begin{center}
 $\Omega:$   $1=P_0<P_1<\cdots<P_i<\cdots<P_{n-1}< P_n=P$
\end{center}
of $P$ is called $c$-normal maximal chain in $G$ if  $P_i$  is $c$-normal in $G$ for each $i\in\{1,2,\ldots,n-1\}$.
 \end{defi}
 
In this paper,  we obtain some new characterizations of $p$-nilpotent groups and supersolvable groups by assuming that Sylow subgroup of $G$ possesses $c$-normal maximal chain in $G$. Our main results are as follows.
 
\begin{thm}
 Let $G$ be a finite group and $p$   the smallest prime divisor of $|G|$.  Suppose that  $P$ is  a Sylow $p$-subgroup of $G$.   If $P$ possesses a $c$-normal maximal chain in $G$, then $G$ is $p$-nilpotent.
\end{thm}

\begin{thm}
Let $G$ be a finite group and $p$ be the smallest prime divisor of $|G|$. Suppose that $N$ is a normal subgroup of $G$  such that $G/N$ is $p$-nilpotent.  If a Sylow $p$-subgroup $P$  of $N$
  possesses a $c$-normal maximal  chain in $G$,  then $G$ is $p$-nilpotent.
\end{thm}

\begin{thm}
 Let $G$ be a finite  group. Suppose that every Sylow subgroup $P$ of $G$ has a $c$-normal maximal 
chain in $G$.  Then $G$ is supersolvable.
\end{thm}

All unexplained notations and terminologies
    are standard and can be found in \cite {DO, HB}.

\section{Proof  of the main results}

To prove the main results, we first need the following lemma.

 \begin{lem} \cite  {WANG}
Let $G$ be a finite group. Suppose that $X \leq H \leq G$ and $N\unlhd G$.

(a) If $X$ is $c$-normal in $G$, then $X$ is also $c$-normal in $H$.

(b) If $X$ is $c$-normal in $G$, then $XN/N$ is $c$-normal in $G/N$.
\end{lem}

Furthermore, we also need the following lemma.

\begin{lem}
Let $G$ be a finite group and $P$ be a Sylow $p$-subgroup of $G$ for some $p\in\pi(G)$. Suppose $N$ is a normal subgroup of $G$.  If  $P$ has a $c$-normal maxima chain in $G$, then  $PN/N$    has a $c$-normal  maximal chain in $G/N$.
\end{lem} 

{\noindent\bf  Proof.} It is trivial when $P\leq N$. So we can assume that $P\not\leq N$. Let 
\begin{center}
$\Omega:$ $1=P_0<P_1<\cdots<P_i<\cdots<P_{n-1}< P_n=P$ 
\end{center}
 be a  maximal chain of $P$ which is $c$-normal  in $G$. By Definition 1.1,  $P_i$ is $c$-normal in $G$ for each $i\in\{1,2,\ldots,n-1\}$. So $P_iN/N$ is $c$-normal in $G/N$ by Lemma 2.1(b). Let us investigate the following subgroup series of $PN/N$ as follows:
\begin{center}
$\overline{\Omega}:$ $N/N= P_0N/N \leq P_1N/N\leq\cdots\leq P_iN/N \leq\cdots\leq P_{n-1}N/N\leq P_nN/N = PN/N$.
\end{center}
For each $i = 1,\ldots , n$, we see that either $P_{i-1}N/N = P_iN/N$ or $P_{i-1}N/N$ is  a maximal subgroup of $P_iN/N$. Therefore, after removing
the equal terms in the above series, we obtain a   maximal chain of $PN/N$ which is $c$-normal in $G/N$. $\Box$\\
 
We now proceed to the proof of Theorem 1.2.\\

{\noindent\bf  Proof of Theorem 1.2.}  Assume that the theorem is false and let $G$ be a counterexample of minimal order. Our proof will be divided into the following six steps. \\

Step 1. $|P|\geq p^2$.

If $|P|\leq p$, then $P$ is cyclic. By Burnside theorem, $G$ is $p$-nilpotent as $p$ is the smallest prime dividing the order of $G$, a contradiction. So $|P|\geq p^2$, and therefore $P$ has non-trivial maximal chains. \\

Step 2. For any normal subgroup $N\neq 1$ of $G$, $G/N$ is $p$-nilpotent.

Consider the quotient group $G/N$.  If $G/N=1$, the assertion is trivial.  Thus we may assume that $G/N\neq1$.  Clearly, $PN/N\in Syl_p(\overline{G})$.   By Lemma 2.2,  $PN/N$ has a $c$-normal maximal chain in $G/N$,  and so $G/N$ satisfies the theorem' hypothesis. By choice of $G$, we conclude that $G/N$ is $p$-nilpotent.\\

Step 3. $G$ has a unique minimal normal subgroup.

Let $N_1$ and $N_2$ be two distinct minimal normal subgroups. Then $N_1\cap N_2=1$. By Step 2, both $G/N_1$ and $G/N_2$ are $p$-nilpotent.  Consequently,
\begin{center}
$G=G/(N_1\cap N_2)\lessapprox G/N_1\times G/N_2$
\end{center}
 and hence $G$ is  $p$-nilpotent, a contradiction. Hence $G$ has a unique minimal normal subgroup. \\

Step 4. $O_{p'}(G)=1$.

 Suppose that $O_{p'}(G)\neq 1$. By Step 2, $G/O_{p'}(G)$ is $p$-nilpotent and hence $G$ is $p$-nilpotent, a contradiction. Therefore,  $O_{p'}(G)=1$.  \\

Step 5. Let $N$ be a unique minimal normal subgroup. Then $N\leq P_{n-1}$.

For convenience, we write $H=P_{n-1}$. By theorem' hypothesis, $H=P_{n-1}$ is $c$-normal in $G$. Therefore, there exists a normal subgroup $K$ of $G$ such that $G=HK$ and $H\cap K\leq H_G$.  We claim that $H_G\neq 1$. Otherwise, suppose that $H_G=1$. Then $H\cap K=1$. Consequently,  
\begin{center}
$|K|_p=|K:K\cap H|_p=|KH:H|_p=|G:H|_p=p$.
\end{center}
By Burnside theorem, $K$ is $p$-nilpotent. Let $M$ be a normal $p$-complement of $K$. Then $M$ is also a normal $p$-complement of $G$ and so $G$ is $p$-nilpotent, a contradiction. Therefore, $H_G\neq1$. Since $N$ be the unique minimal normal subgroup of $G$, we have $N\leq H_G\leq H=P_{n-1}$. \\

Step 6. $|N|=p$, the final contradiction.

By Step 5, we know that $N\leq P_{n-1}$. So there exists a integer $k\in\{1,2,\ldots,n-1\}$ such that $N\leq P_k$, but $N\not\leq P_{k-1}$. Since $N$ is the unique minimal normal subgroup of $G$, we have $P_{{k-1}_G}=1$ (otherwise, $N\leq P_{{k-1}_G}\leq P_{k-1}$, a contradiction).
If $k=1$, then $|N|=|P_1|=p$, and the conclusion holds.  Suppose that $k\geq 2$. Observe that $|P_k:P_{k-1}|=p$, we have $P_k=P_{k-1}N$. Consequently,
\begin{center}
$|N:N\cap P_{k-1}|= |NP_{k-1}:P_{k-1}|=|P_k:P_{k-1}|=p$. 
\end{center}
Furthermore, by theorem hypothesis, $P_{k-1}$ is $c$-normal in $G$, so there exist a normal subgroup $M$ of $G$ such that $G=P_{k-1}M$ and $P_{k-1}\cap M\leq P_{{k-1}_G}=1$. Since $N$ is the unique minimal normal subgroup of $G$,  $N\leq M$ and hence  $N\cap P_{k-1}=1$. This implies that $|N|=|N:N\cap P_{k-1}|=p$.

 Moreover, by $N/C$-theorem, $G/C_G(N)\lessapprox Aut(N)$.  Hence $|G/C_G(N)|\bigm | p-1$. Since $p$ is the smallest prime dividing the order of $G$, we have $G=C_G(N)$ and hence $N\leq Z(G)$. So $G/Z(G)$ is $p$-nilpotent  and hence  $G$ is $p$-nilpotent. This is a final contradiction.
 The proof of theorem is complete.
$\Box$\\

Next, we apply Theorem 1.2 to prove Theorem 1.3.\\

{\noindent\bf  Proof of Theorem 1.3.} 
By the theorem' hypothesis, $G/N$ is $p$-nilpotent. Let $T/N$   be a normal $p$-complement of $G/N$. Since
$Syl_p(N)=Syl_p(T)$,  it follows that $P\in Syl_p(T)$. 
  The theorem' hypothesis also provides a maximal chain of subgroups of $P$ as follows:
  \begin{center}
$\Omega:$ $1=P_0<P_1<\cdots<P_i<\cdots<P_{n-1}< P_n=P$ 
\end{center}
such that $P_i$ is $c$-normal in $G$ for each $i\in\{1,2,\ldots,n-1\}$.   By Lemma 2.1(a), every
$P_i$   is $c$-normal in $T$ for $i\in\{1,2,\ldots,n-1\}$.
  Applying Theorem 1.2, we conclude that $T$ is $p$-nilpotent. 
 Let $T_0$  be     the normal $p$-complement of $T$. Then 
 $T_0$ is a Hall $p'$-subgroup of $G$.
  Moreover, $T_0$  is characteristic in $T$, and $T$ is normal in $G$. 
This implies  that $T_0$ is normal in $G$.  Consequently, $G$ is $p$-nilpotent.  The proof of theorem is complete.
 $\Box$\\
 
Finally, we prove Theorem 1.4 by induction on the order of $G$.\\

{\noindent\bf  Proof of Theorem 1.4.} Assume that the theorem is false and let $G$ be a counterexample of minimal order.
Let $p$ be the smallest prime divisor of $|G|$. Then $G$ is $p$-nilpotent by Theorem 1.2. Let $U$ be a Hall normal
$p$-subgroup of $G$. It is easy to see  that $U$ satisfies the theorem hypothesis by Lemma 2.1 (a). By induction,
$U$ is supersolvable, and hence $G$ possesses Sylow tower property of supersolvable type. Let $q$ be the
largest prime divisor of $|G|$ and $Q$ is a Sylow $q$-subgroup of $G$. Then $Q$ is normal in $G$. Applying
Lemma 2.2, $G/Q$ satisfies the theorem hypothesis. By induction, $G/Q$ is supersolvable.

Let $N$ be a minimal normal subgroup of $G$ contained in $Q$. By Lemma 2.2, $G/N$ satisfies the theorem hypothesis, and hence $G/N$ is supersolvable  by  the minimality of $G$. Since the class of all supersolvable groups is a saturated formation, $N$ is the unique minimal normal subgroup of $G$ contained in $Q$.

It remains to show  that $|N|=q$. 
By theorem hypothesis, $Q$ has a maximal chain 
\begin{center}
$1=Q_0<Q_1<\cdots<Q_i<\cdots<Q_{m-1}<Q_m=Q$.
\end{center}
which is $c$-normal in $G$. 

First, suppose that $N\not\leq Q_{m-1}$. Since $N$ is the unique minimal normal subgroup of $G$ contained in $Q$, we have $Q_{{m-1}_G}=1$. 
Moreover, since  $Q_{m-1}$ is $c$-normal in $G$,   there exist a normal subgroup $M$ of $G$ such that $G=Q_{m-1}M$ and $Q_{m-1}\cap M\leq Q_{{m-1}_G}=1$. It follows that $Q\cap M$ is also normal $q$-subgroup of $G$ and $Q\cap M\in Syl_q(M)$. Hence $N\leq Q\cap M$. Furthermore,  
\begin{center}
$|M|_q=|M:M\cap Q_{m-1}|_q=|MQ_{m-1}:Q_{m-1}|_q=|G:Q_{m-1}|_q=|Q:Q_{m-1}|=q$.
\end{center}
 This implies that $|N |=|Q\cap M|=|M|_q=q$.  
 
Second, suppose that $N\leq Q_{m-1}$. 
Let $k$ be the smallest number with $N\leq Q_k$, that is,   $N\leq Q_k$, but $N\not\leq Q_{k-1}$. Since $N$ is the unique minimal normal subgroup of $G$ contained in $Q$, we have $Q_{{k-1}_G}=1$ (otherwise, $N\leq Q_{{k-1}_G}\leq Q_{k-1}$, a contradiction). If $k=1$, then $|N|=|Q_1|=q$, as desired.  Thus  we may assume that $k\geq 2$. Since $|Q_k:Q_{k-1}|=q$, we have $Q_k=Q_{k-1}N$. Consequently,
\begin{center}
$|N:N\cap Q_{k-1}|= |NQ_{k-1}:Q_{k-1}|=|Q_k:Q_{k-1}|=q$.
 \end{center}
 Furthermore, by theorem hypothesis, $Q_{k-1}$ is $c$-normal in $G$, so there exist a normal subgroup $M$ of $G$ such that $G=Q_{k-1}M$ and $Q_{k-1}\cap M\leq P_{{k-1}_G}=1$. Since $N$ is the unique minimal normal subgroup of $G$ contained in $Q$, we have  $N\leq Q\cap M$, and hence  $N\cap Q_{k-1}=1$. This implies that $|N|=|N:N\cap P_{k-1}|=q$.

Finally, since $q$ is the largest prime divisor of $|G|$ and $G/N$ is supersolvable, we  conclude that $G$ is supersolvable. 
This is final contradiction.  The proof of theorem is complete. $\Box$





\bibliographystyle{amsplain}

\end{document}